\documentclass[12pt,reqno]{amsart}

\usepackage[arrow,matrix,curve]{xy}

\usepackage[dvips]{graphicx} 

\usepackage{amssymb, latexsym, amsmath, amscd, array, hyperref
}

\newtheorem{theorem}{Theorem}[section]

\newtheorem{fact}[theorem]{Fact}

\theoremstyle{definition}

\numberwithin{equation}{section}
\newcommand\N {{\mathbb N}} 

\newcommand\R {{\mathbb R}}
\newcommand\Q {{\mathbb Q}}

\newcommand\astr{{}^\ast\R}

\newcommand\Los{{\L}o{\'s}}

\newcommand\GOT{{\mbox{\rm GOT\,$\setminus$PATHOS}}}

\newcommand\IST{{\mbox{\rm IST}}}

\author[AG]{Alexander E. Gutman}\address{A. Gutman, Sobolev Institute
of Mathematics, Siberian Branch Russian Academy of Sciences,
Novosibirsk State University, Novosibirsk, Russia}
\email{gutman@math.nsc.ru}

\author[MK]{Mikhail G. Katz}\address{M. Katz, Department of
Mathematics, Bar Ilan University, Ramat Gan 52900
Israel}\email{katzmik@macs.biu.ac.il}

\author[TK]{Taras S. Kudryk}\address{T. Kudryk, Department of
Mathematics, Lviv National University, Lviv, Ukraine}
\email{kudryk@mail.lviv.ua}

\author[SK]{Semen S. Kutateladze}\address{S. Kutateladze, Sobolev
Institute of Mathematics, Novosibirsk State University, Russia}
\email{sskut@math.nsc.ru}

\newbox\GrossOneBox
\newbox\grossOneBox
\setbox\GrossOneBox\hbox{\raise-.4pt\hbox{${\rm O}\mskip-10.2mu\raise-2.7pt\hbox{$^1$}\mskip3mu$}}
\setbox\grossOneBox\hbox{\raise-4.5pt\hbox{\small$\mskip-1.3mu{\hbox{\rm\footnotesize O}}\mskip-10.0mu\raise-3.0pt\hbox{$^1$}\mskip1mu$}}
\def\GrossOne{{\copy\GrossOneBox}}

\begin{document}

\thispagestyle{empty}

%\huge

\title{\emph{The Mathematical Intelligencer} flunks the Olympics}

%Keywords: Archimedean axiom; hyperreals; incomparable quantities;
%infinity; Robinson.

\maketitle

\begin{abstract}
\emph{The Mathematical Intelligencer} recently published a note by
Y. Sergeyev that challenges both \emph{mathematics} and
\emph{intelligence}.  We examine Sergeyev's claims concerning his
purported \emph{Infinity computer}.  We compare his \emph{grossone}
system with the classical Levi-Civita fields and with the hyperreal
framework of A.~Robinson, and analyze the related algorithmic issues
inevitably arising in any genuine computer implementation.  We show
that Sergeyev's \emph{grossone} system is unnecessary and vague, and
that whatever consistent subsystem could be salvaged is subsumed
entirely within a stronger and clearer system (IST).  Lou Kauffman,
who published an article on a \emph{grossone}, places it squarely
\emph{outside} the historical panorama of ideas dealing with infinity
and infinitesimals.
\end{abstract}

\tableofcontents

\section{Grossone olympics}

In the summer of 2015, some of us were approached by an editor of
\emph{The Mathematical Intelligencer (TMI)} with a request to respond
to a piece of what they felt was \emph{pseudo-science}, published
without their knowledge in \emph{TMI}.\, As noted in
\cite[p.~393]{Da15}, I.~Grattan-Guinness argued that ``the demarcation
between science and pseudo-science is not clearly drawn.''  While
agreeing with Grattan-Guinness, in the present article we argue that
in some cases the demarcation is drawn clearer than in others.

Yaroslav Sergeyev has developed a positional system for infinite
numbers in numerous articles over the past decade.  By 2015,
\emph{MathSciNet} listed 19 such articles, starting with \cite{Se03}.
His ``Olympic Medal'' note \cite{Se15a} in \emph{TMI} purports to be
an application of his \emph{grossone} system to ranking countries
lexicographically according to the number of gold, silver, and bronze
medals they earned in the olympics.  Sergeyev's system is closely
related to the field of rational functions in one variable and to the
classical Levi-Civita field, with a non-Archimedean structure provided
by a suitable lexicographic ordering (a more detailed comparison with
the Levi-Civita fields appears in Section~\ref{s43}).

Sergeyev appears to be making claims of significant progress in the
field of nonstandard models.  The reaction of the experts to
Sergeyev's claims has been lukewarm.  Joel David Hamkins, a leading
authority on mathematical logic and foundations, reacted as follows to
Sergeyev's claims: ``It seems to me that there is very little that is
new in this topic, and basically nothing to support the grand claims
being made about it.''  \cite{Ha15} In this text, we will analyze
Sergeyev's claims in more detail.

Shamseddine's group has used Levi-Civita fields to develop computer
implementations exploiting infinite numbers (see Section~\ref{s42}),
without engaging in the sort of rhetorical \emph{flou artistique} that
envelopes a typical Sergeyev performance.  Pure and applied
mathematicians may sometimes use different standards of rigor but
Sergeyev's case is a rather different problem.

Nonstandard models of arithmetic were developed as early as 1933 by
Skolem using purely constructive methods (in particular not relying on
any version of the axiom of choice); see e.g., \cite{Sk33},
\cite{Sk34}, \cite{Sk55}, and \cite[Section~3.2]{KKM}.

Conservative extensions of the Peano axioms (PA) were studied in
\cite{Kr69} and \cite{HKK}.  Subsequently \cite{HK86} described both a
family of nonstandard versions of PA itself, and $n$-th order PA for
different values of $n$, that are conservative extensions of PA itself
and respectively $n$-th order PA (see Proposition 2.3 there), and also
nonstandard versions containing additional stronger saturation axioms,
that are not conservative extensions (see Theorem 3.2 there).  All of
these theories are conservative with respect to ZFC, as is IST (see
Section~\ref{s6}).

\cite{Av05} showed how to use weak theories of nonstandard arithmetic
to treat fragments of calculus and analysis.  If (as apparently
claimed in \cite{Lo15}) what Sergeyev is attempting to do is develop
such nonstandard models, he is certainly doing it without
acknowledging prior work in the field.

Contrary to Sergeyev's earlier announcements, Nobel Prize laureate
Robert Aumann will \emph{not} be attending Sergeyev's june '16 meeting
in Italy.

\section{Transfering the sine function}
\label{s1}

A few years ago, one of the authors asked Sergeyev through email what
the sine of his \emph{grossone} was, and he replied that it is
\[
\sin(grossone).
\]
The author in question did not have the heart to ask Sergeyev what
\[
\sin^2(grossone) + \cos^2(grossone)
\]
is, and how exactly his ``infinity computer" can know it other than
being told case-by-case about every possible identity in mathematics.
The point is that neither the field of rational functions nor
Sergeyev's \emph{grossone} system possesses a transfer principle (see
below) or any equivalent procedure.

In his list of areas where his ideas are claimed to be potentially
fruitful, Sergeyev mentions differential equations.  Surely for this
he will need to know that the sine function is defined on the extended
system with its usual properties.  This is what makes the question
about~$\sin(grossone)$ crucial.

The \emph{transfer principle} is a type of theorem that, depending on
the context, asserts that rules, laws or procedures valid for a
certain number system, still apply (i.e., are ``transfered'') to an
extended number system.  Thus, the familiar extension~$\Q\subseteq\R$
preserves the property of being an ordered field.  To give a negative
example, the extension
\[
\R\subseteq\R\cup\{\pm\infty\}
\]
of the real numbers to the so-called \emph{extended reals} does not
preserve such a property.  The hyperreal extension
\[
\R\subseteq\astr
\]
preserves \emph{all} first-order properties, including the
trigonometric identity~$\sin^2 x + \cos^2 x =1$ (valid for all
hyperreal~$x$, including infinitesimal and infinite values
of~$x\in\astr$).  For a more detailed discussion, see the textbook
\emph{Elementary Calculus} \cite{Ke86}.

The revolutionary idea that there does exist a system, sometimes
called \emph{hyperreal numbers}, satisfying such a transfer principle
is due to the combined effort of \cite{He48}, \cite{Lo55}, and
\cite{Ro61}, and has roots in Leibniz's \emph{Law of continuity} and
his distinction between \emph{assignable} and \emph{inassignable}
numbers; see \cite{KS2}, \cite{KS1}, \cite{Ba16a}, \cite{Ba16b}, as
well as \cite{Bl16a}.  We will provide an explanation of the
extension~$\R\subseteq\astr$ in Section~\ref{s5}.

Sergeyev sometimes grudgingly acknowledges the debt to Robinson.
However, in many publications Sergeyev unfortunately presents the idea
as his own, as noted by Vladik Kreinovich in his \emph{MathSciNet}
review of Sergeyev's book \cite{Kr03}.  Peter W. Day's review of
Sergeyev's article at \cite{Da06} mentions the connection to the
transfer principle, lacking in Sergeyev's system.  Additional critical
reviews are \cite{Zl09} and \cite{Ku11}.

Sergeyev himself introduces his symbol for infinity in the following
terms:
\begin{quote}
A new infinite unit of measure has been introduced for this purpose as
the number of elements of the set~$\N$ of natural numbers.  It is
expressed by the numeral~\GrossOne{} called \emph{grossone}.  It is
necessary to note immediately that~\GrossOne{} is neither
Cantor's~$\aleph_0$ nor~$\omega$.  Particularly, it has both cardinal
and ordinal properties as usual finite natural numbers
\cite[p.~8101]{Se12}.
\end{quote}
It is easy to detect serious logical problems with such a definition.
Sergeyev's claim that his \GrossOne{} has both cardinal and ordinal
properties is a purely declamative pronouncement.  A reader might have
expected such a claim in a refereed mathematical periodical to be
justified by a clever definition, but it is not.  As it stands,
Sergeyev's claim is merely a thinly veiled admission of an
inconsistency, couched in an attempt to dress up a bug to look like a
feature.  Similarly, Sergeyev's attempted definition of \GrossOne{} as
somehow ``the number of elements of the set~$\N$'' contradicts other
passages where \GrossOne{} is \emph{included} as a member of~$\N$,
resulting in an embarrassing circularity.%
\footnote{\label{f1}See further on circularity of Sergeyev's
definitions in footnote~\ref{f5}.}

The point we wish to emphasize is that the plausibility that such a
scheme might actually work after being sufficiently cleaned-up of
superfluous pathos%
\footnote{\label{f1b}The English word \emph{pathos} is etymologically
related to~$\pi\acute{\alpha}\theta o\varsigma$, passion.}
(including inconsistencies), is entirely due to Robinson's insights
implementing Leibniz's ideas about the distinction between assignable
and inassignable numbers, on the one hand, and implementing Leibniz's
law of continuity as the transfer principle, on the other.

In his writings, Sergeyev introduces his \emph{grossone}, announces
that it is infinite, and blithely assumes that anything algebraic, or
even from analysis, that can be done with ordinary numbers can be done
when the \emph{grossone} is adjoined.  Such mathematical assertions
require proof, which are lacking in the analyzed note.

\section{Debt to Robinson}
\label{s2}

The tendency to give insufficient credit to Robinson is clearly on
display in the ``Olympic medal'' as the reference to Robinson's theory
is concealed in an obscure phrase in such a way that an uninformed
reader will be unable to gauge its significance.

For the benefit of such a reader, we provide the following
clarification.  As far as providing a lexicographic ordering for the
olympic medals are concerned, it would be sufficient to take the
\emph{grossone} to be equal to a number~$p$ greater than the total of
all the medals attributed at the olympics, for example~$p$ equal a
million, and work with number representation in base~$p$.  Then
obviously~$p$ will satisfy all the usual rules governing finite
numbers, because~$p$ itself is a finite number.  However, Sergeyev's
system is obviously not tailor-made for the games.  Rather, the
alleged significance of Sergeyev's system is its purported
applicability to a broad range of scientific problems, without any
apriori limitation on the size of the sample.  For this reason he
wishes to use an infinite \emph{grossone} value for~$p$.  In fact, the
ordinary rational numbers suffice for this purpose, as we explain in
Section~\ref{s7}.

This is where his (pseudo)mathematical claims become questionable.
His framework presupposes a number system which \emph{properly}
extends the usual one, yet obeys the usual laws, i.e., a transfer
principle (see Section~\ref{s1}).  But Sergeyev's system does not obey
a transfer principle in any mathematically identifiable form, as
Sergeyev appears to acknowledge in his~$\sin(grossone)$ comment.  The
\emph{grossone} calculator will be able to compute values necessary
for scientific work only to the extent that one or another version of
the transfer principle is successfully implemented.  While Robinson's
system does obey a transfer principle, Sergeyev is sparing in
acknowledging his debt to Robinson.

Thus, in his keynote address in Las Vegas '15, Sergeyev declares that
\begin{quote}
The new computational \emph{methodology is not related to} the
non-standard analysis and gives the possibility to execute
computations of a new type simplifying fields of Mathematics where the
usage of infinity and/or infinitesimals is required.  \cite{Se15b}
(emphasis added)
\end{quote}
This strikes us as a somewhat economical way of acknowledging
intellectual indebtedness.  It is as if someone proclaimed himself to
be the inventor of relativity theory and declared that his
``methodology is not related to" the work of Albert Einstein.

Sergeyev's infringement on Robinson's framework appears to be
tolerated by the decision-makers in the mathematics community, in a
way that would not be tolerated if the infringement were in a field
like differential geometry or Lie theory.  An infringement upon
Robinson's framework is tolerated at least in part because the field
created by Robinson has been marginalized, not least through the
(combined) efforts of Paul Halmos and Errett Bishop (see e.g.,
\cite{KK11}, \cite{KK12a}, \cite{Ka15}), and of Connes (see
\cite{KKM}, \cite{KL}).  As a result, a number of Robinson's students
were unable to obtain positions at PhD-granting institutions in the
1970s.  An additional factor seems to be Robinson's apparent
insistence that logic has to take a more prominent place in graduate
programs in mathematics, provoking animosity on the part of some
mathematicians.

Robinson's framework is a fruitful modern research area that has
attracted many researchers.  Thus, Terry Tao developed certain
arguments on approximate groups exploiting ultraproducts that would be
difficult to paraphrase without them.  The ultraproducts form a bridge
between discrete and continuous analysis, and enable a unified
framework for a treatment of both Hilbert's fifth probem and Gromov's
theorem on groups of polynomial growth; see \cite{Ta14} for details.

\section{Comparison with work by other scholars}

In this section we will compare Sergeyev's work with that of other
scholars, in chronological order.

\subsection{Levi-Civita fields}
\label{s43}

David Tall used Levi-Civita fields under the name \emph{superreal} to
popularize teaching calculus via infinitesimals in \cite{Ta79}.
Levi-Civita fields is a classical topic with a long history.  It was
studied in \cite{RL}.  Sergeyev exploits his \emph{grossone} in place
of the variable~$x$ in the Levi-Civita fields with the lexicographic
ordering, but comments that
\begin{quote}
Levi-Civita numbers are built using a generic infinitesimal
$\varepsilon$ \ldots whereas our \emph{numerical} computations with
[in]finite quantities are concrete and not generic.
\cite[p.~2]{Se15c} (emphasis added)
\end{quote}
Two years earlier, Sergeyev compared the concrete \emph{grossone}
numeral to Levi-Civita in the following terms (we make no attempt to
correct the grammar):
\begin{quote}
${}^5$ At the first glance the \emph{numerals} (7) can remind numbers
from the Levi-Civita field (see [20]) that is a very interesting and
important precedent of algebraic manipulations with infinities and
infinitesimals. However, the two mathematical objects have several
crucial differences. They have been introduced for different purposes
by using two mathematical languages having different accuracies and on
the basis of different methodological foundations. In fact,
Levi-Civita does not discuss the distinction between numbers and
\emph{numerals}.  His numbers have neither cardinal nor ordinal
properties; they are build [sic] using a generic infinitesimal and
only its rational powers are allowed; he uses symbol~$\infty$ in his
construction; there is no any \emph{numeral} system that would allow
one to assign \emph{numerical} values to these numbers; it is not
explained how it would be possible to pass from \ldots a generic
infinitesimal~$h$ to a concrete one (see also the discussion above on
the distinction between numbers and \emph{numerals}). In no way the
said above should be considered as a criticism with respect to results
of Levi-Civita. The above discussion has been introduced in this text
just to underline that we are in front of two different mathematical
tools that should be used in different mathematical contexts.
\cite[p.~10671, note~5]{Se13} (emphasis added)
\end{quote}
Sergeyev's use of the terms \emph{numeral} (both as adjective and
noun) and \emph{numerical} is vague.  Certainly real numbers cannot be
used in computer implementations, and one needs to work instead with a
specific representation such as decimals.  Shamseddine and his
colleagues are surely aware of this in their work with the Levi-Civita
fields (see Section~\ref{s42}).

Sergeyev has a talent for turning \emph{pathos}%
\footnote{See our etymological comment in footnote~\ref{f1b}.}
into patent.  Affected pathos was also characteristic of the
\emph{superior ideology} of the former Soviet Union where he was
raised.  Sergeyev seems to have learned the lesson of the rhetorical
effectiveness of \emph{superior ideology}.  Levi-Civita may have done
the same mathematics a hundred years earlier than Sergeyev, but the
former says a mere ``$x$'' and the latter says a superior
``\emph{numeral},'' ergo the latter is on so much higher an
ideological plane.

\subsection{Shamseddine's work on Levi-Civita fields}
\label{s42}

A group of researchers around K. Shamseddine have been developing
software based on the Levi-Civita field for handling certain
calculations with infinity and infinitesimals; see e.g., the article
\cite{Sh15} and \url{http://www.bt.pa.msu.edu/index_cosy.htm}

These scholars typically refrain from assorting their work with the
kind of rhetoric that typically accompanies a Sergeyev performance,
such as:
\begin{enumerate}
\item
Sergeyev does not acknowledge properly indebtedness to Robinson,
particularly in the matter of the transfer principle (see
Section~\ref{s1}), painting himself as a pioneer in the area.
\item
Sergeyev does not acknowledge properly that what he is working with is
a version of the classical Levi-Civita fields, seeking to emphasize
what he claims to be the novelty of his system.
\item
Sergeyev seeks to spice up his writing with an assortment of colorful
\emph{principles} that have little bearing on an actual computer
implementation, such as his stylized insistence on the part being less
than the whole.
\end{enumerate}
With regard to this last point, \cite{BD} developed a mathematical
theory of \emph{numerosities} to express this idea mathematically, but
its Sergeyevan incarnation seems to have little mathematical content.

\subsection{Kauffman on \GrossOne}

L.~Kauffman is a leading topologist today.  The \emph{Kauffman
bracket} \cite{KL94} is a staple of 3-manifold invariants.  His
article ``Infinite computations and the generic finite" \cite{Kau15}
uses Sergeyev's notation \GrossOne.  Sergeyev managed to cite this
recent paper of Kauffman's already in three texts.  Thus, Sergeyev
sends the reader to Kauffman (and other texts) ``In order to see the
place of the new approach in the historical panorama of ideas dealing
with infinite and infinitesimal'' \cite[p.~24]{Se16}.  However,
Kauffman himself clearly distances himself from Sergeyev's
``methodology" in the following terms:
\begin{quote}
In my paper about the Grossone, I point out that the logic of this
formalism is identical (in my version) to using~$1+x+x^2+\ldots+x^G$
as a \emph{finite sum} with~$G$ a generic positive integer.  One can
then manipulate the series and look at the limiting behaviour in many
cases. There is no need to invoke any new concepts about infinity.
This point of view may be at variance with the interpretations of
Yaroslav [Sergeyev] for his invention, but I suggest that this is what
is happening here.  \cite{Kau15b}
\end{quote}
In no way can Kauffman's work or comments be interpreted as support
for Sergeyev.  Nor does Kauffman place Sergeyev ``in the historical
panorama'' etc., contrary to Sergeyev's claim.  Quite the opposite,
Kauffman writes that ``[t]here is no need to invoke any new concepts
about infinity,'' thereby placing Sergeyev squarely \emph{outside} a
``historical panorama of ideas dealing with the infinite.''

\section{The hyperreal extension}
\label{s5}

In an approach to analysis within Robinson's framework, one works
with the pair~$\R\subseteq\astr$ where~$\R$ is the usual ordered
complete \emph{Archimedean} continuum, whereas~$\astr$ is a proper
extension thereof.  A proper extension of the real numbers could be
called a \emph{Bernoullian} continuum, in honor of Johann Bernoulli
who was the first systematically to use an infinitesimal-enriched
continuum as the foundation for analysis.  For historical background
see \cite{BK}, \cite{Ba13}, \cite{Ba14}, \cite{Ka15b}.  The
extension~$\astr$ obeys the transfer principle (see Section~\ref{s1}).

The field~$\astr$ is constructed from~$\R$ using sequences of real
numbers.  The main idea is to represent an infinitesimal by a sequence
tending to zero.  One can get something in this direction without
reliance on any nonconstructive foundational material.  Namely, one
takes the ring of all sequences, and quotient it by the equivalence
relation that declares two sequences to be equivalent if they differ
only on a finite set of indices.

The resulting object is a proper ring extension of~$\R$, where~$\R$ is
embedded by means of the constant sequences.  However, this object is
not a field.  For example, it has zero divisors.  But if one quotients
it further in such a way as to obtain a field (by extending the kernel
to a \emph{maximal} ideal), then the quotient will be a field, called
a hyperreal field.

To motivate the construction further, it is helfpul to analyze first
the construction of~$\R$ itself using sequences of rational numbers.
Let~$\Q^\N_C$ denote the ring of Cauchy sequences of rational numbers.
Then
\begin{equation}
\label{51}
\R=\Q^\N_C/\text{MAX}
\end{equation}
where ``MAX'' is the maximal ideal in~$\Q^\N_C$ consisting of all null
sequences (i.e., sequences tending to zero).

The construction of a Bernoullian field can be viewed as refining the
construction of the reals via Cauchy sequences of rationals.  This can
be motivated by a discussion of rates of convergence as follows.  In
the above construction, a real number~$u$ is represented by a Cauchy
sequence~$\langle u_n : n\in\N\rangle$ of rationals.  But the passage
from~$\langle u_n\rangle$ to~$u$ in this construction sacrifices too
much information.  We seek to retain some of the information about the
sequence, such as its ``speed of convergence."  This is what one means
by ``relaxing" or ``refining" the equivalence relation in the
construction of the reals from sequences of rationals.

When such an additional piece of information is retained, two
different sequences, say~$\langle u_n\rangle$ and~$\langle
u'_n\rangle$, may both converge to~$u\in\R$, but at different speeds.
The corresponding ``numbers" will differ from~$u$ by distinct
infinitesimals.  If~$\langle u_n\rangle$ converges to~$u$ faster than
$\langle u'_n\rangle$, then the corresponding infinitesimal will be
smaller.  The retaining of such additional information allows one to
distinguish between the equivalence class of~$\langle u_n\rangle$ and
that of~$\langle u'_n\rangle$ and therefore obtain distinct hyperreals
infinitely close to~$u$.  For example, the sequence
$\langle\frac{1}{n^2}\rangle$ generates a smaller infinitesimal than
$\langle\frac{1}{n}\rangle$.

A formal implementation of the ideas outlined above is as follows.
Let us present a construction of a hyperreal field~$\astr$.
Let~$\R^{\N}$ denote the ring of sequences of real numbers, with
arithmetic operations defined termwise.  Then we have
\begin{equation}
\label{52}
\astr=\R^{\N}\!/\text{MAX}
\end{equation}
where ``MAX'' is a suitable maximal ideal.  What we wish to emphasize
is the formal analogy between \eqref{51} and \eqref{52}.  In both
cases, the subfield is embedded in the superfield by means of constant
sequences.

We now describe a construction of such a maximal ideal exploiting a
suitable finitely additive measure~$m$.  The ideal MAX consists of all
``negligible'' sequences~$\langle u_n\rangle$, i.e., sequences which
vanish for a set of indices of full measure~$m$, namely,
\[
m\big(\{n\in\N:u_n=0\}\big)=1.
\]
Here~$m:\mathcal{P}(\N)\to\{0,1\}$ (thus~$m$ takes only two
values,~$0$ and~$1$) is a finitely additive measure taking the
value~$1$ on each cofinite set,%
\footnote{For each pair of complementary \emph{infinite} subsets
of~$\N$, such a measure~$m$ ``decides'' in a coherent way which one is
``negligible'' (i.e., of measure~$0$) and which is ``dominant''
(measure~$1$).}
where~$\mathcal{P}(\N)$ is the set of subsets of~$\N$.  The
subset~$\mathcal{F}_m\subseteq\mathcal{P}(\N)$ consisting of sets of
full measure~$m$ is called a free ultrafilter.  These originate with
\cite{Ta30}.  The construction of a Bernoullian continuum outlined
above was therefore not available prior to that date.

The construction outlined above is known as an ultrapower
construction.  The first construction of this type appeared in
\cite{He48}, as did the term \emph{hyper-real}.  The transfer
principle (see Section~\ref{s1}) for this extension is an immediate
consequence of the theorem of \Los; see \cite{Lo55}.

\section{A detailed technical report on GOT}
\label{s6}

The analysis presented in this section is an extension of the report
\cite{GK}.  We formulate our analysis in the framework of Nelson's
Internal Set Theory (IST) first presented in \cite{Ne77}.

The difference between Nelson's approach and Robinson's can be
illustrated in the context of the underlying number system as follows.
Robinson extended the real number field to a hyperreal number field
with infinitesimals (for example, by the ultrapower approach of
Section~\ref{s5}).  In contrast with Robinson's approach, Nelson
proceeded axiomatically and revealed both infinitesimals and illimited
numbers within the real number field itself.%
\footnote{\label{f5}This point seems to have escaped Sergeyev, who
claims it to be an advantage of the grossone system that the infinite
numbers are found \emph{within}~$\N$, allegedly unlike nonstandard
analysis; see \cite[p.~95, note~3]{CD}.  Elsewhere Sergeyev claims
that, on the contrary, \GrossOne{} is ``the number of elements in
$\N$'', leading to a circularity already mentioned in
footnote~\ref{f1}.}
To this end, Nelson introduced a new one-place predicate ``to be
standard" together with the appropriate axioms.  Both Nelson's and
Robinson's theories are conservative extensions of the traditional
foundational framework of the Zermelo--Fraenkel set theory.  For
further discussion see \cite{KK}.

\subsection{Logical status of Sergeyev's theory}

Sergeyev's reasoning is not only informal but often vague and
inaccurate.  The inaccuracies include his definition of the grossone
as ``the number of elements in set of natural numbers'' (which may
appeal to the uneducated but mathematically speaking is nonsensical),
as well as his delphic pronouncements as to ``the whole being greater
than the part" and the distinction between ``numbers and numerals"
(see Section~\ref{s43}).  Such superfluous PATHOS needs to be removed
before a consistent theory can be identified.  A reader with some
mathematical culture can give formal shape to Sergeyev's postulates,
as done in \cite{GK} to some extent.  The result is a formal theory of
signature
\[
S = \big\{=, \in, \GrossOne{} \big\}
\]
(here~$\in$ is the membership relation while \GrossOne{} is the
grossone).  We will abbreviate the theory as \GOT{}.  Here ``GOT''
stands for GrossOne Theory, while ``PATHOS'' alludes to the
inconsistencies of Sergeyev's system and his efforts to sweep them
under the rug by means of \emph{le flou artistique} via affected
pathos or passionate enthusiasm; see Section~\ref{s43}.  Thus, \GOT{}
is the axiomatic formal theory in the language of signature~$S$ whose
axiomatic background is given by all of Sergeyev's postulates, both
explicitly stated and implicitly assumed in his papers.

\begin{fact}
Each axiom of \GOT{} is a trivial consequence of the axioms of any
classical nonstandard set theory, provided~$\GrossOne$ is understood
as the factorial of an infinitely large integer.
\end{fact}

This is shown in \cite{GK}. In particular, the axioms of \GOT{} are
easily proven in Nelson's \IST{}, with \GrossOne{} evaluated as the
factorial of an arbitrary infinitely large natural number.  Therefore,

\begin{fact}
The theory \GOT{} is weaker than \IST{}.
\end{fact}

By definition, the theory is weaker whenever is has fewer theorems.
Note that, for formal theories, \emph{weaker} does not mean
\emph{worse}; nor does \emph{stronger} mean \emph{better}.  For
instance, a theory whose theorems are all the statements, i.e., an
inconsistent theory, is the strongest one, but it is hardly the best
one.  Nevertheless, in certain circumstances, a weaker theory cannot
be regarded as \emph{new} as compared to a stronger theory.

\begin{fact}
\GOT{} is not a new theory.
\end{fact}

Indeed, \GOT{} is weaker than a well-known theory, \IST{}, and
moreover, the axioms of \GOT{} are easily proven in
\IST{}. Consequently, any reasoning within \GOT{} can be automatically
converted into the corresponding and almost identical reasoning in
\IST{}. In particular, \GOT{} cannot prove any \emph{new} result,
since each result proven in \GOT{} is already a result of a well-known
theory.  Actually, even proofs within \GOT{} cannot be \emph{new},
since every such proof is almost identical to an automatically
produced proof in a well-known theory.

\begin{fact}
\GOT{} is dramatically weaker than \IST{}.
\end{fact}

It suffices to note that \IST{} features a powerful and fruitful tool
known as the Transfer Principle (see Section~\ref{s1}), which is
absent from the theory \GOT{}. In addition, \GOT{} has no analogs of
Idealization and Standardization Principles, which makes it almost
impossible to prove any serious assertion in \GOT{} without appealing
to informal or implicit assumptions.

\begin{fact}
Consistency of \GOT{} is not justified by its originator.
\end{fact}

In many of Sergeyev's papers, one cannot find a single attempt
formally to justify the consistency of the grossone theory.  Only due
to \cite{GK} do we know that \GOT{} is consistent relative to \IST{}
(see also \cite{Va12}).  Furthermore, employing the fact that \IST{}
is consistent relative to ZFC (see Nelson's article \cite{Ne77}) and
that ZFC is consistent relative to ZF (a result of Goedel's; see his
constructible universe \cite{Go38}), we may conclude that \GOT{} is
consistent relative to the standard set theory. (This is however not
surprising, since \GOT{} is weaker than a well-known relatively
consistent theory.)

It is good to know which facts a theory can prove, but for a theory to
be useful it also very important to know which facts it cannot prove.
To become a generally accepted legitimate mathematical tool, a theory
should be unable to prove strange or \emph{pathological} results.  The
corresponding formal property of a theory is called
\emph{conservativity}.

By definition, a theory T* of signature S* is a conservative extension
of a weaker theory T with smaller signature S whenever T* has exactly
the same theorems in signature S as T has. Suppose that we have a
generally accepted theory T (say, ZFC) and let a new theory T* (say,
\IST{}) extend T and introduce new primary notions (in our example, the
notion of standard set).  The fact that T* is a conservative extension
of T means the following: if T* allows us to prove some result R and R
does not involve new primary notions, then R is not pathological, as
it can also be proven in the generally accepted theory T.  Therefore,
any conservative extension of a customary theory can be (and should
be) accepted as a legitimate mathematical tool.  Namely, it has the
same deduction strength and every \emph{sensible} fact it can prove
can be proven by usual means, without any new axioms or new notions.

\IST{} is known to be a conservative extension of ZFC, as shown by
Powell's theorem presented in \cite{Ne77}.  This nontrivial and very
important fact makes \IST{} a generally accepted mathematical theory.

\begin{fact}
The question of conservativity of \GOT{} is ignored by its originator.
\end{fact}

Again, only due to \cite{GK} do we know that \GOT{} is weaker than
\IST{}, which, in its turn, is a conservative extension of ZFC. Hence,
so is \GOT{}: if a set-theoretic fact can be proven in \GOT{}, it can
also be proven in \IST{} and, thus, in ZFC. Without knowing this, even
a consistent theory need not be accepted.

Therefore, without employing nontrivial facts from contemporary
nonstandard analysis, Sergeyev's reasoning remains a powerless,
informal, weak theory with doubtful consistency, which cannot be
generally accepted due to its doubtful conservativity. On the other
hand, if we employ the facts from nonstandard analysis, the grossone
theory turns out to be merely a powerless and weak theory which cannot
be regarded as new.

\subsection{Algorithmic status of Sergeyev's theory}
\label{s62}

An algorithmic problem is the task of finding an algorithm which,
given a constructive object as input, produces a constructive object
as output so that the output is related to the input in a desired way,
and this fact is provable within a suitable theory under
consideration. Therefore, solvability and complexity of an algorithmic
problem depends on the underlying theory.

A solution to an algorithmic problem is an algorithm supplied with a 
justification, i.e., with a proof (within a theory) of the assertion that 
the algorithm works correctly and actually solves the problem. On the other 
hand, a weaker theory has fewer proofs (which is a direct consequence of 
the definition) and thus fewer solvable algorithmic problems.

\begin{fact}
Within a weaker theory, there are more unprovable and undecidable
statements, more unsolvable algorithmic problems, while solutions to
solvable algorithmic problems are more complex.
\end{fact}

Recalling that \GOT{} is weaker than \IST{}, we conclude the following.

\begin{fact}
\label{f8}
Each algorithmic problem unsolvable in \IST{} is similarly unsolvable
in \GOT{}; if an algorithmic problem has a complex solution in \IST{},
it either has an even more complex solution in the system \GOT{} or is
even unsolvable in \GOT{}.
\end{fact}

Furthermore, being a conservative extension of ZFC, \IST{} has exactly
the same solvable set-theoretical problems as ZFC has.  This
circumstance allows us to derive the following fact.

\begin{fact}
\label{f9}
Every unsolvable set-theoretical problem is unsolvable in \GOT{};
solvable set-theoretical problems are more complex or even unsolvable
in \GOT{}.
\end{fact}

There is a number of problems listed in \cite{GK} which encounter
certain theoretical obstacles to finding an algorithmic solution.
Some of the problems are \IST{}-specific, other are purely
set-theoretical or analytical.  According to facts \ref{f8} and
\ref{f9} we have the following fact.

\begin{fact}
Each of the algorithmic problems enumerated in the article \cite{GK}
is either more complex or even unsolvable in \GOT{}.
\end{fact}

\subsection{Specific algorithmic problems concerning grossone}

Within \GOT{}, the main tool is the ``positional system with
base~\GrossOne{}'' in which the role of numerals is played by
``multilevel polynomials" in a single variable denoted \GrossOne{},
with rational coefficients and exponents. We will refer to these
polynomials as \emph{grossnumerals}.  They are multilevel in the sense
that the exponents (power indices) need not be numbers and may also be
(multilevel) polynomials.  Every grossnumeral has finite height.
Suitable formal definitions are presented in \cite{GK} (and are absent
from Sergeyev's papers).

If we restrict the height of grossnumerals to 1, we obtain the usual
polynomials in one variable. The algorithmic problems in the classical
calculus of such polynomials are far from being new. They are all
solved, long ago and completely. Anything new can occur only under
consideration of numerals having arbitrary finite height.

The set of grossnumerals cannot be called a ``calculus" unless it is
supplied with a set of algorithms which implement such key operations
as reduction to canonical form and comparison. Without such
algorithms, one cannot speak of any computer realization of the
calculus, either.

The important point here is that the implementation of the basic
calculus operations in the set of grossnumerals encounters certain
theoretical obstacles in \IST{} and ZFC.  According to Section~\ref{s62},
they encounter even more serious problems in the weaker \GOT{}.  The
issues are thoroughly described in \cite{GK}, and the main problem is
as follows.

\begin{fact}
There is no known algorithm that, given grossnumerals x and y, would
determine which of the following holds true:~$x<y$,~$x=y$, or~$x>y$.
\end{fact}

The latter problem must be solved in order to be able to speak of a
\emph{calculus}, for otherwise we would not be able to perform such
elementary procedures as reducing similar terms or listing the terms
in descending order by their degree.  Nevertheless, algorithmic
solvability of these procedures remains unknown. The corresponding
hypothesis is based on rather nontrivial facts on o-minimality and
decidability of the order structure of reals with exponent (see
bibliographic references [11] and [13] in \cite{GK}).

Thus, currently there is no algorithm able to compare grossnumerals
or, for that matter, to check the inequalities
\[
1<\GrossOne^{{\GrossOne}^{-1}} < 2.
\]
Such an algorihm could hardly appear in any of Sergeyev's papers.
Indeed, he provides the following characterisation of infinite
numbers: ``Infinite numbers in this numeral system are expressed by
numerals having at least one grosspower grater [sic] than zero.''
\cite[p.~60]{Se07} But the grossexponent~${\GrossOne}^{-1}$ is indeed
greater than zero; yet the number~$\GrossOne^{{\GrossOne}^{-1}}$ must
be infinitely close to~$1$ if even a most rudimentary form of the
transfer principle (see Section~\ref{s1}) is to be satisfied.  Yet
according to Sergeyev's characterisation,
$\GrossOne^{{\GrossOne}^{-1}}$ would turn out to be ``infinite".
Whenever Sergeyev's assertions are specific enough to be checked, one
finds errors, including freshman calculus level errors.  

This particular error appeared in ``Blinking fractals'' \cite{Se07}
published in \emph{Chaos, Solitons, and Fractals}, and was
subsequently criticized in \cite{GK}.  Sergeyev blinked and modified
his text in a number of online databases, so as to remove the error,
including its current ResearchGate version.  As of 2015, no official
correction whatsoever appeared in \emph{Chaos, Solitons, and
Fractals}.

This episode indicates how far removed the questions under
consideration are from any computer implementation.  The comparison
problem is completely ignored in Sergeyev's papers, and this is not
surprising: the problem is challenging even in \IST{}, while in \GOT{}
it is much more complex due to the absence of a suitable transfer
principle.

With the above taken into account, it becomes clear why all
screenshots of a calculator presented in Sergeyev's papers contain
only grossnumerals of height 1.

\begin{fact}
An actual grossone calculator does not exist.
\end{fact}

Grossnumerals of height 1 are just ordinary polynomials of one
variable, and software for the corresponding calculus is commonplace
nowadays.  Contemporary symbolic computation packages provide much
more sophisticated machinery.  The grossone theory is so poorly
designed and underdeveloped that a toy calculator is the only tool
which can be created on its basis.

\section{Olympic ranks need no ``numerical infinities''}
\label{s7}

In his note ``The Olympic medals, ranks, lexicographic ordering, and
numerical infinities,'' Sergeyev represents the basics of grossone
theory (as he does in each of his numerous papers containing the
symbol~$\GrossOne$) under the pretext of applying it to a
``mathematical problem'' related to the lexicographic ranking method.
The problem is caused by the fact that, contrary to other known
ranking methods, the lexicographic method does not assign numerical
ranks to various medal distributions, it only orders them, i.e.,
determines which distribution is higher and which is lower.  Sergeyev
suggests using grossnumerals as ``numerical'' ranks of arbitrary medal
distributions and emphasizes that his suggestion solves the problem
without upper bounds on the number of medals awarded by a single
country as well as on the number of the medal classes (gold, silver,
etc.).

We will demonstrate that the approach suggested by Sergeyev is useless
and any application of a theory of infinite numbers is overkill for
such a trivial aim.  Indeed, the lexicographic order can be made
numerical in a very easy, reasonable, and practical way by means of
ordinary standard rational numbers.

Suppose that there are infinitely (but countably) many medal classes.
List them in descending order and associate with successive natural numbers:
$1$~for~``gold,\!''~$2$~for~``silver,\!''~$3$~for~``bronze,\!''
$4$,~$5$,~$6$,~etc.~for all the rest.
Each competitor can win an arbitrary finite set of medals
which can be encoded by a~finite word with positive integers as ``letters.\!''
For instance, the word~$w=\langle5,0,12,1\rangle$ encodes the fact that a~competitor has won
$5$~medals of class~$1$,~$0$~medals of class~$2$,~$12$~medals of class~$3$,~$1$~medal of class~$4$, and~$0$~medals of any other class.
The task is to invent a practical method (an algorithm) of calculating a number~$R(w)$ for any word~$w$ in such a way that
the equality~$R(u)>R(v)$ be equivalent to~$u\succ v$, where~$\succ$ is the lexicographic order on words:
$$
\text{$u\succ v\ \Leftrightarrow\ u_1=v_1$,~$\dots$,~$u_{n-1}=v_{n-1}$,~$u_n>v_n$ for some~$n$.}
$$
(Here~$w_n$ is the~$n\mskip2mu$th letter of a word~$w$, with~$w_n=0$ for~$n$ greater than the length of~$w$.)

The method proposed by Sergeyev consists in defining the ``numerical'' rank~$R_S(w)$
of a~word~$w=\langle w_1,\dots,w_L\rangle$ of length~$L$ as the grossnumeral
\[
R_S(w)\ =\ w_1\GrossOne^{L-1} + w_2\GrossOne^{L-2} + \cdots +
w_{L-1}\GrossOne^{1} + w_L\GrossOne^{0}.
\]
How useful is such a solution, however?  Sergeyev regards~$R_S(w)$ as
a ``numerical'' rank just because it is a ``number'' in the sense of
his grossone theory.  Both theoretically and practically, this is
nothing but a mere replacement of a word~$\langle
w_1,\dots,w_L\rangle$ with a more bulky expression of the form
$w_1\GrossOne^{L-1} + \cdots + w_L\GrossOne^{0}$.  This expression
cannot be written in any other numerical form and cannot be used in
any software other than the hypothetical ``Infinity Calculator'' based
on the mythical ``Infinity Computer technology.\!''

We will now indicate a very simple and honest method of solving the
above-stated ``problem.\!''  Note first that, for the aim under
consideration, there is no need for any artificial numbers, and the
standard rational numbers with their standard order are undoubtedly
sufficient.  This is so because, as is well known, every countable
linear order embeds into the standard ordered set of rationals, and
this is true, in particular, for the lexicographically ordered set of
words which represent medal distributions.  So, the task is merely in
choosing a specific order-preserving rational encoding of the words.
The encoding can be as simple as follows.  Given a word~$w=\langle
w_1,\dots,w_L\rangle$, set
\[
 R(w)\ =\ \sum_{n=1}^L
 2^{-(w_1+\cdots+w_{n-1}+n-1)}\sum_{m=1}^{w_n}2^{-m}.
\]
Here~$R(w)\in[0,1)$ is the rational number whose binary representation
(representation in the positional numeral system with base~$2$) has
the form
$$
  0{\mskip2mu\boldsymbol.}
      \underset{w_1\text{ ones}}{\underbrace{\strut 11...1}} 0
      \underset{w_2\text{ ones}}{\underbrace{\strut 11...1}} 0\ \ldots\ 0
      \underset{w_L\text{ ones}}{\underbrace{\strut 11...1}}\,.
$$

It is an easy exercise to show that the encoding~$R$ meets the
required condition, i.e., assigns greater ranks~$R(w)$ to
lexicographically greater words~$w$.  Note also that medal
distributions are uniquely (and easily) determined by their numerical
ranks.  It is also worth observing that~$R$~reflects certain emotional
aspects related to medals wins: the awarding of the first medal of
a~given class is felt as a more exciting and significant achievement
than awarding the second one, and so on. This circumstance results in
the fact that the medal distributions with close numerical ranks are
also ``psychologically'' close.

As an illustration, we present the 2014 Winter Olympics medal table of
competitors (in lexicographic order) and their medal distributions
supplemented with the corresponding exact binary ranks, and
approximate decimal ranks.

\bigskip
\centerline{\textit{2014 Winter Olympics medal table}}
\medskip

{\def\m#1{\hbox to .8ex{\hss#1\hss}}
\begin{tabular}{|l|ccc|l|l|}\hline{\vline height 2.5ex width 0pt}%
 {\bf Country} & \multicolumn{3}{|c|}{\bf Medals} & {\bf Binary} & {\bf Decimal} \\\hline{\vline height 3ex width 0pt}%
 Russia        & \ \m{13} & \m{11} &  \m{9} \ & 0.11111111111110111111111110111111111 & 0.9999389 \\
 Norway        & \ \m{11} &  \m{5} & \m{10} \ & 0.1111111111101111101111111111        & 0.9997520 \\
 Canada        & \ \m{10} & \m{10} &  \m{5} \ & 0.111111111101111111111011111         & 0.9995114 \\
 United States & \  \m{9} &  \m{7} & \m{12} \ & 0.111111111011111110111111111111      & 0.9990196 \\
 Netherlands   & \  \m{8} &  \m{7} &  \m{9} \ & 0.11111111011111110111111111          & 0.9980392 \\
 Germany       & \  \m{8} &  \m{6} &  \m{5} \ & 0.111111110111111011111               & 0.9980311 \\
 Switzerland   & \  \m{6} &  \m{3} &  \m{2} \ & 0.1111110111011                       & 0.9915771 \\
 Belarus       & \  \m{5} &  \m{0} &  \m{1} \ & 0.11111001                            & 0.9726562 \\
 Austria       & \  \m{4} &  \m{8} &  \m{5} \ & 0.1111011111111011111                 & 0.9686870 \\
 France        & \  \m{4} &  \m{4} &  \m{7} \ & 0.11110111101111111                   & 0.9677658 \\
 Poland        & \  \m{4} &  \m{1} &  \m{1} \ & 0.11110101                            & 0.9570312 \\
 China         & \  \m{3} &  \m{4} &  \m{2} \ & 0.11101111011                         & 0.9350585 \\
 South Korea   & \  \m{3} &  \m{3} &  \m{2} \ & 0.1110111011                          & 0.9326171 \\
 Sweden        & \  \m{2} &  \m{7} &  \m{6} \ & 0.11011111110111111                   & 0.8745040 \\
 Czech Republic& \  \m{2} &  \m{4} &  \m{2} \ & 0.1101111011                          & 0.8701171 \\
 Slovenia      & \  \m{2} &  \m{2} &  \m{4} \ & 0.1101101111                          & 0.8583984 \\
 Japan         & \  \m{1} &  \m{4} &  \m{3} \ & 0.1011110111                          & 0.7412109 \\
 Finland       & \  \m{1} &  \m{3} &  \m{1} \ & 0.1011101                             & 0.7265625 \\
 Great Britain & \  \m{1} &  \m{1} &  \m{2} \ & 0.101011                              & 0.6718750 \\
 Ukraine       & \  \m{1} &  \m{0} &  \m{1} \ & 0.1001                                & 0.5625000 \\
 Slovakia      & \  \m{1} &  \m{0} &  \m{0} \ & 0.1                                   & 0.5000000 \\
 Italy         & \  \m{0} &  \m{2} &  \m{6} \ & 0.0110111111                          & 0.4365234 \\
 Latvia        & \  \m{0} &  \m{2} &  \m{2} \ & 0.011011                              & 0.4218750 \\
 Australia     & \  \m{0} &  \m{2} &  \m{1} \ & 0.01101                               & 0.4062500 \\
 Croatia       & \  \m{0} &  \m{1} &  \m{0} \ & 0.01                                  & 0.2500000 \\
 Kazakhstan    & \  \m{0} &  \m{0} &  \m{1} \ & 0.001                                 & 0.1250000 \\[1ex]\hline
\end{tabular}}

\bigskip\bigskip

\section{Publication venue}
\label{s8}

This rebuttal did not appear in the journal \emph{The Mathematical
Intelligencer} where Sergeyev's note originally appeared because five
successive versions of our rebuttal were rejected by that journal, in
spite of at least one favorable referee report.

\section{Conclusion}

The Olympic medals ranking was considered in Sergeyev's note in
\emph{The Mathematical Intelligencer} without any serious mathematical
treatment.  The note's shortcomings include serious issues of
attribution of prior work.

\section*{Acknowledgments}

We are grateful to Rob Ely for helpful suggestions.  We thank the
anonymous referee for \emph{Foundations of Science} for helpful
comments.  M.~Katz was partially funded by the Israel Science
Foundation grant no.~1517/12.

\medskip\noindent \textbf{Alexander E. Gutman} was born in 1966 in
Novokuznetsk, USSR.  He is the Head of the Laboratory of Functional
Analysis at the Sobolev Institute of Mathematics in Novosibirsk and a
professor at Novosibirsk State University.  He authored four books
and over 70 papers in functional analysis and Boolean valued analysis.

\medskip\noindent \textbf{Mikhail G. Katz} (B.A. Harvard University,
'80; Ph.D. Columbia University, '84) is Professor of Mathematics at
Bar Ilan University.  His book \emph{Systolic geometry and topology}
was published by the American Mathematical Society.  In his spare time
he is more likely than not to be defending the small, the tiny, and
the infinitesimal.

\medskip\noindent \textbf{Taras S. Kudryk} (born 1961, Lviv, Ukraine)
is a Ukrainian mathematician and associate professor of mathematics at
Ivan Franko National University of Lviv. His main interests are
nonstandard analysis and its applications to functional analysis. He
is the author of books about nonstandard analysis (in Ukrainian and
English) and textbooks about functional analysis (in Ukrainian)
co-authored with Vladyslav Lyantse.  Kudryk has performed research in
nonstandard analysis in collaboration with V.~Lyantse and V\'\i tor
Neves.  His publications appeared in \emph{Matematychni Studii},
\emph{Siberian Journal of Mathematics}, and \emph{Logica Universalis}.

\medskip\noindent
\textbf{Semen S. Kutateladze} was born in 1945 in Leningrad (now
St.~Petersburg).  He is a senior principal officer of the Sobolev
Institute of Mathematics in Novosibirsk and professor at Novosibirsk
State University.  He authored more than 20 books and 200 papers in
functional analysis, convex geometry, optimization, and nonstandard
and Boolean valued analysis.  He is a member of the editorial boards
of \emph{Siberian Mathematical Journal}, \emph{Journal of Applied and
Industrial Mathematics}, \emph{Positivity}, \emph{Mathematical Notes},
etc.
\end{document}